\documentclass[12pt,openbib]{article}
\usepackage{amsmath}
\usepackage[cp1251]{inputenc}
\usepackage[russian, english]{babel}
\usepackage{amsfonts}
\usepackage{amsfonts,amssymb}
\usepackage{amssymb}
\usepackage{latexsym}
\usepackage{euscript}
\usepackage{enumerate}
\usepackage{graphics}
\usepackage[dvips]{graphicx}

\newtheorem{teo}{\bf Theorem}

\newtheorem{state}{Statement}
\newtheorem{opr}{Definition}

\newtheorem{task}{Problem}

\begin{document}

\begin{center}
{\bf Parities on 2-knots and 2-links \\
D.\,A.~Fedoseev\footnote[1]{%
{\it Denis Fedoseev} --- V. A. Trapeznikov Institute of Control Sciences of RAS,
e-mail: denfedex@yandex.ru.},
V.\,O.~Manturov\footnote[2]{%
{\it Vassily Manturov} --- Chelyabinsk State University,
e-mail: vomanturov@yandex.ru}.
}

\vspace{3mm}

{
2-dimensional knots and links are studied in the article. The notion of parity is introduced via techniques similar to the ones used by the second named author in 1-dimensional case. By using parity new invariants are constructed and known invariants are refined.} 

\end{center}

\section{Introduction}

The foundation for {\em parity theory} was laid in the paper \cite{knot_parity} for one-dimensional (virtual) knots. The main idea behind it is: all crossings of a diagram of a (virtual) knot can be naturally divided into {\em even} ones and {\em odd} ones so that this parity behaves well under Reidemeister moves. It is shown that it is possible to refine many knot invariants via parity as well as to create new invariants valued in knot diagrams, see \cite{imn}. Due to existance of ``picture-valued'' (that is, diagram-valued) invariants the following principle holds for virtual knots: \\ 

{\it if a diagram is complicated enough, it can be found as a subdiagram in any equivalent diagram}. \\

Virtual knots admit many different parities. The simpliest of them --- the Gaussian --- comes from Gauss (chord) diagram of a knot: every crossing corresponds to a chord of the Gauss diagram; a crossing is called {\it even} if the corresponding chord intersects an even number of chords, and {\it odd} in other case. If the given knot is classical, then all chords of its Gauss diagram are even. This fact was known to Gauss himself. Therefore, the Gaussian parity on classical knots is trivial. Moreover it is proved that all parities are trivial on classical knots. For virtual knots the Gaussian parity is defined exactly the same way and is non-trivial.

The Gaussian parity is not the only one existing in virtual knot theory; it is possible to define a set of {\em parity axioms} so that the theorems and invariants obtained thanks to parity hold for any parity (that is a map from the set of crossings to the set $\{$even, odd$\}$) satisfiyng the axioms.

Non-trivial parities exist not only on virtual knots and allow one to prove theorems and create invariants in different theories. For example consider the so-called {\it free knots}. In the classical knot theory every crossing is decorated with ``over/under'' structure telling which arc of the crossing is over (its coordinate along the projection direction greater) and which is under. Thanks to the parity even without this structure one can create strong invariants. Free knots are a rough factorisation of virtual knots obtained by ``forgetting'' the ``over/under'' structure in the crossings. Initially Turaev \cite{tur} suggested a conjecture that free knots are trivial. By using partiy arguments, the second named author disproved this conjecture in \cite{imn}.

There are many different approaches leading to different parities. One of the natural sources of parities is (co-)homology of the ambient space. Informally one can say that in the classical case there are no non-trivial parities since plane has trivial homology groups. Nevertheless, many parity theoretical methods can be applied to the classical case thanks to a number of topological constructions (see \cite{chri-man, kras-man, gnk}). \\

Parity theory for 1-knots (classical, virtual, etc.) boils down to the decoration of crossings, that is codimension 1 singularities of the projection. A natural idea is to perform the same for higher dimensions. In the present work we study two-dimensional knots. Once again we give parity to codimension 1 singularities --- double lines.

Similarly to classical 1-knots are embeddings of a circle into a 3-sphere, classical 2-knots are embeddings of a 2-sphere into $\mathbb{R}^4$ (or a sphere $S^4$), see Definition \ref{2knot} later. The standard way of describing 2-knots is to consider their diagrams --- two-dimensional complexes of a particular kind. A diagram is a projection in general position of a knot in 4-space onto a 3-subspace where for each double line it is said which leaf is ``over'' and which is ``under'' and for every triple point the same is said for all three leaves incident to that point.

To diagrams are equivalent if they correspond to the same knot; in other words their preimages are isotopic. In one-dimensional case two diagrams represent the same knot if and only if they can be connected by a sequence of Reidemeister moves and trivial isotopies (that is isotopies preserving the combinatorial type of the object). In two-dimensional case it was shown (see \cite{ros}) that two diagrams correspond to the same knot if and only if they can be connected by a sequence of moves from a certain finite set ({\it Roseman moves}). In the classical case it is imposed that all the complexes can be embedded into a 3-space and the moves to be compatible with the embedding.

Relaxing the embedding condition one obtains {\it abstract 2-knots} (see for example \cite{abstract, winter}). Forgetting the ``over/under'' structure one gets {\it free} 2-knots. They turn out not only to be important as a tool for the study of 2-knots in $\mathbb{R}^4$ but are interesting all by themselves.

We will defined parity for abstract 2-knots and 2-links using an analog to the Gauss diagram for 1-knots (see Definition \ref{def:gauss}). Moreover we formulate the {\it 2-parity axioms} based on {\it local links} appearing in Roseman moves. After that the notion of parity is used to refine certain invariants of abstract 2-knots (say, quandle). \\

The structure of the present paper is as follows. In the second Section we give the main definitions of 2-knots and links, their diagrams and moves on them. The third Section is devoted to the definition of parities on 2-knots and the study of their properties. That is in Subsection 3.1 we consider the Gaussian parity in detail; Subsection 3.2 contains a general approach to the definition of parity. A {\it weak parity} is defined in the Subsection 3.3. Subsection 3.4 deals with parities on two-component links. Subsection 3.5 is devoted to projection and hierarchy of parities. In Subsection 3.6 we state the important {\it succession principle} and prove a theorem regarding the succession of equality for certain classes of 2-knots. In Section 4 we consider some invariants of 2-knots, in particular refined via parity. In Section 5 several open problems of 2-knot theory are given. \\

The authors are grateful to S. Carter and I.M. Nikonov for various useful discussions. The authors are especially grateful to D.P. Ilyutko.

This research was supported by a grant from Russian scientific fund (project 16-11-10291).

\section{Basic notions and definitions}

In the present paper we study two-dimensional knots and links. Also we consider {\it abstract} 2-knots. Let us recall the basic definitions of those objects.

\begin{opr}
\label{2knot}
A {\rm 2-knot} (resp. {\rm an $n$-component 2-link}) is a smooth embedding in general position of a 2-sphere $S^2$ (resp. disjoint union on $n$ spheres) into $\mathbb{R}^4$ or $S^4$ up to isotopy.
\end{opr}

Like the usual 1-knots, 2-knots are studied via thier {\it diagrams}. There are several approaches to knot diagrams. We will use the ``2-diagrams in 3-space'' approach. Let us define it following Roseman \cite{ros, ros2} and Winter \cite{winter}. Note that in the same way one can define knot diagrams for any dimension.

Consider a 2-link $L\subset F^3\times [0,1]$ (in the classical case the manifold $F$ equals either $\mathbb{R}^3$ or $S^3$) and a projection $\pi: L\to F$.

The {\it crossing set} $D^*$ of a link $L$ is a closure of the set $\{x\in \pi(L)|\pi^{-1}(x)\ge 2\}$. The {\it double point set} $D$ is defined as $\pi^{-1}(D^*)$. The {\it branch point set} $B$ is a subset of the link $L$ where $\pi$ fails to be an embedding. The set of {\it actual double points} $D_0$ is a subset of $D$ consisting of such $x\in L$ that $|\pi^{-1}(\pi(x))|=2$. We also define the {\it overcrossing set} $D_+$ as the set of $x\in D_0$ such that for any two points $x,y \in \pi^{-1}(\pi(x))$ holds $\pi_I(y)<\pi_I(x)$ with $\pi_I$ denoting the natural projection $F\times [0,1]\to [0,1]$. The {\it undercrossing set} is defined as $D_-=D_0\setminus D_+$.

We will say that a link $L$ is in a {\it general position} with respect to the projection $\pi$ if

\begin{enumerate}
    \item $D$ is a union of embedded closed 1-submanifolds of the link $L$ with normal crossings;
    \item $B$ is a finite subset of the set $D$ and for every $b\in B$ there exists an open neighborhood $U\subseteq D$ such that $U\setminus B$ consists of two components $U_0, U_1$ which are emebdded into $F$ via $\pi$ so that $\pi(U_0)=\pi(U_1)$.
\end{enumerate}

For every link there is an isotopy sending it into a link in general position with respect to $\pi$.

\begin{opr}
Diagram $d$ of a link $L$ in general position with respect to the projection $\pi$ is an image of the link $L$ under the projection $\pi(L)$; every point of the crossing set with exactly two preimages is marked with information which preimage has greater coordinate in $[0,1]$.
\end{opr}

\begin{figure}[h]
\begin{center}
$\includegraphics[width=90mm]{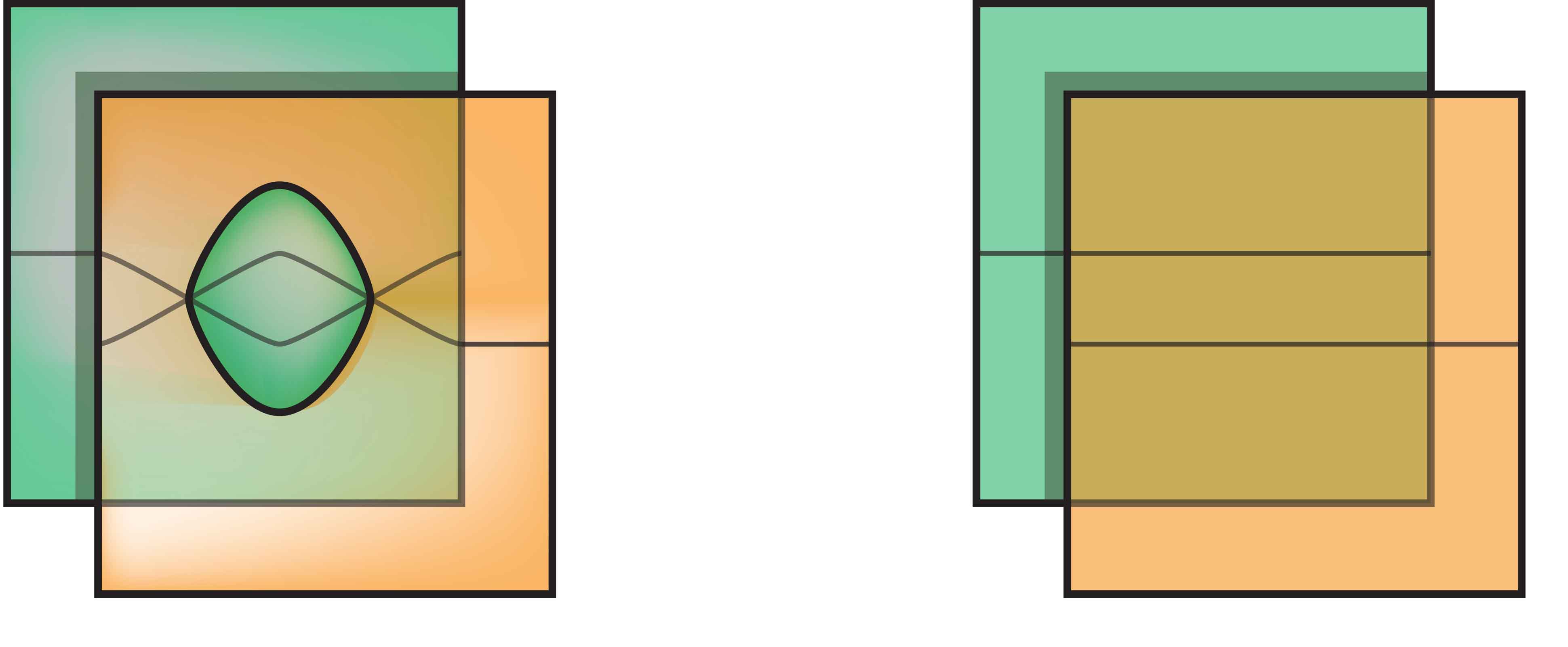}$ 
\caption{First Roseman move $\mathcal{R}_1$: elliptic move of the type $\Omega_2$.} \label{ris:R1}
\end{center}
\end{figure}

\begin{figure}[h]
\begin{center}
$\includegraphics[width=90mm]{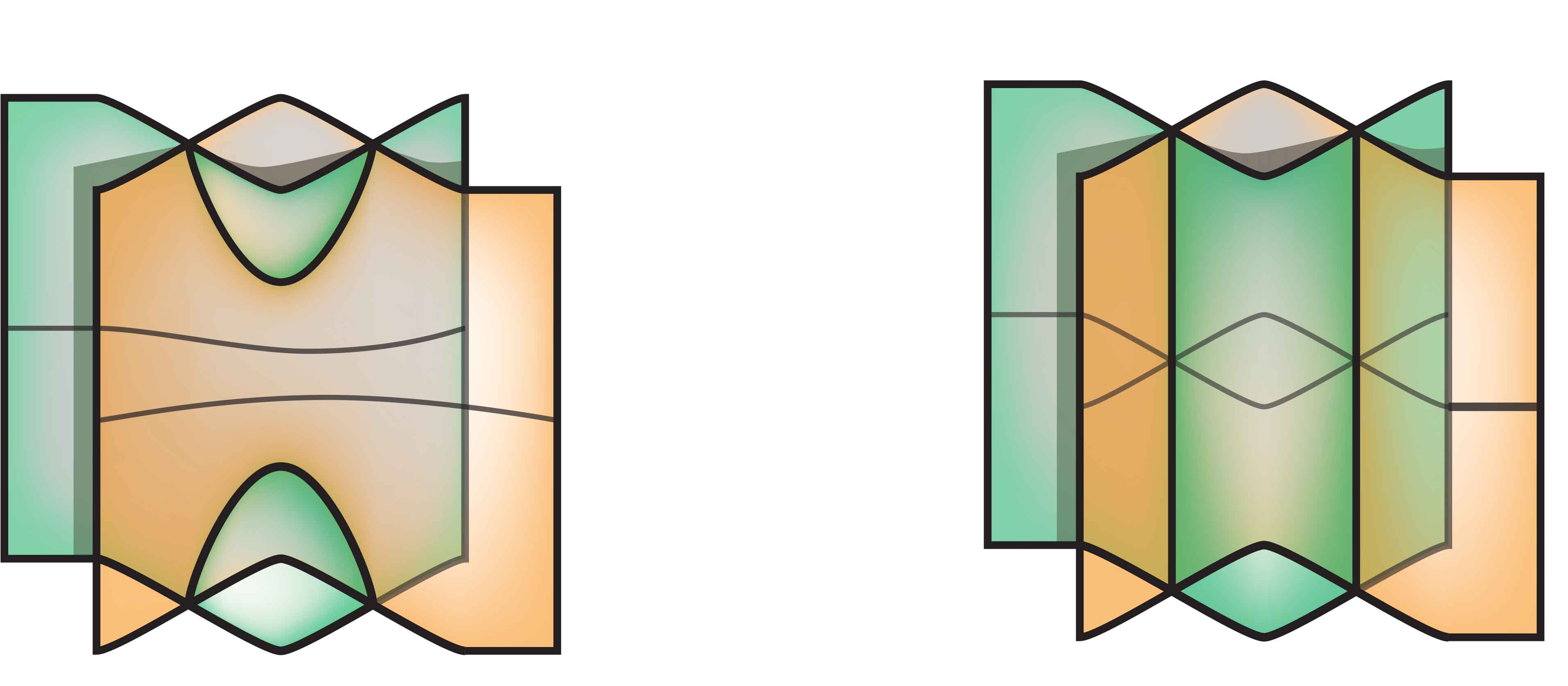}$ 
\caption{Second Roseman move $\mathcal{R}_2$: hyperbolic move of the type $\Omega_2$.} \label{ris:R2}
\end{center}
\end{figure}

\begin{figure}[h]
\begin{center}
$\includegraphics[width=90mm]{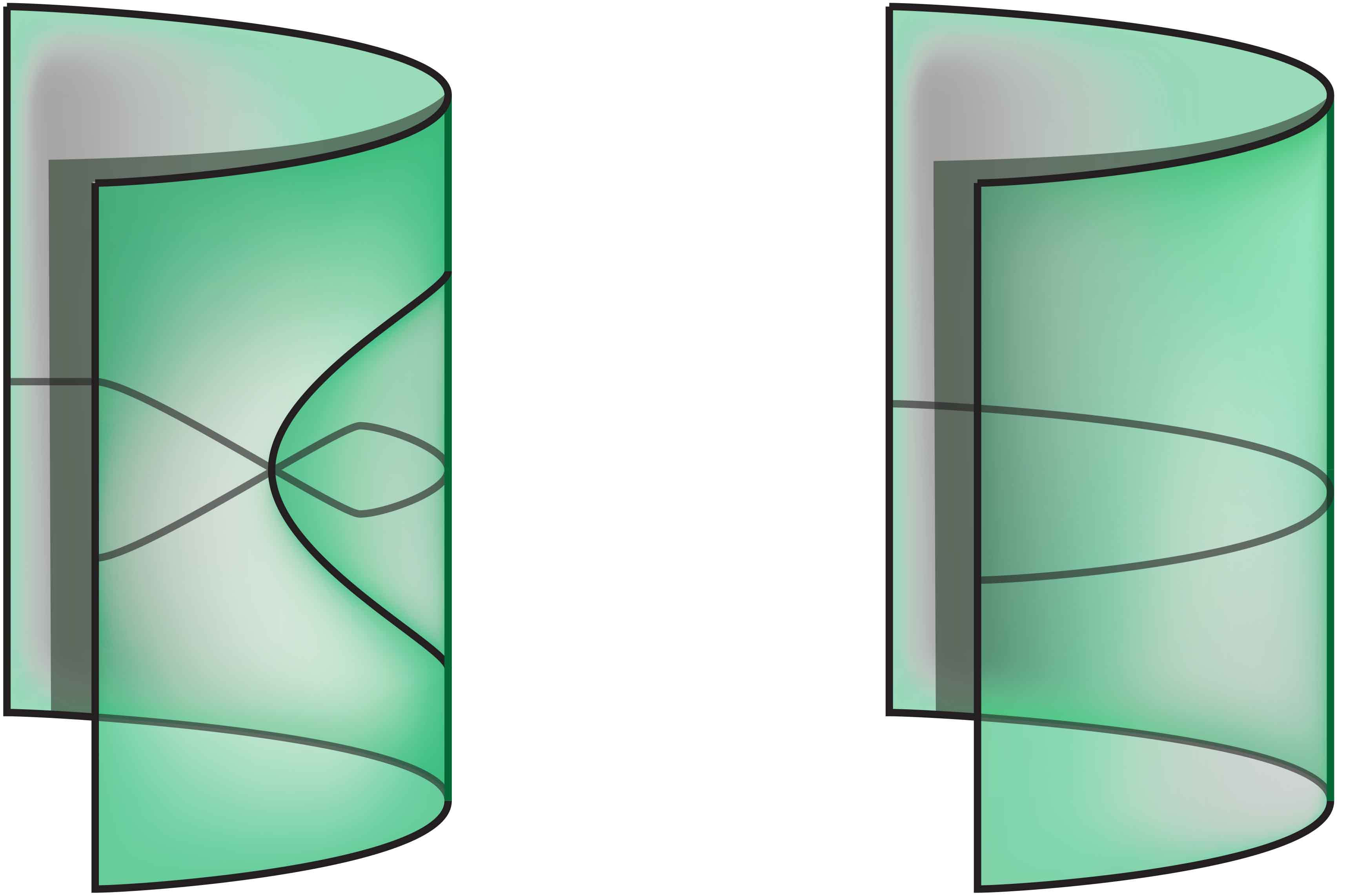}$ 
\caption{Third Roseman move $\mathcal{R}_3$: elliptic move of the type $\Omega_1$.} \label{ris:R3}
\end{center}
\end{figure}

\begin{figure}[h]
\begin{center}
$\includegraphics[width=90mm]{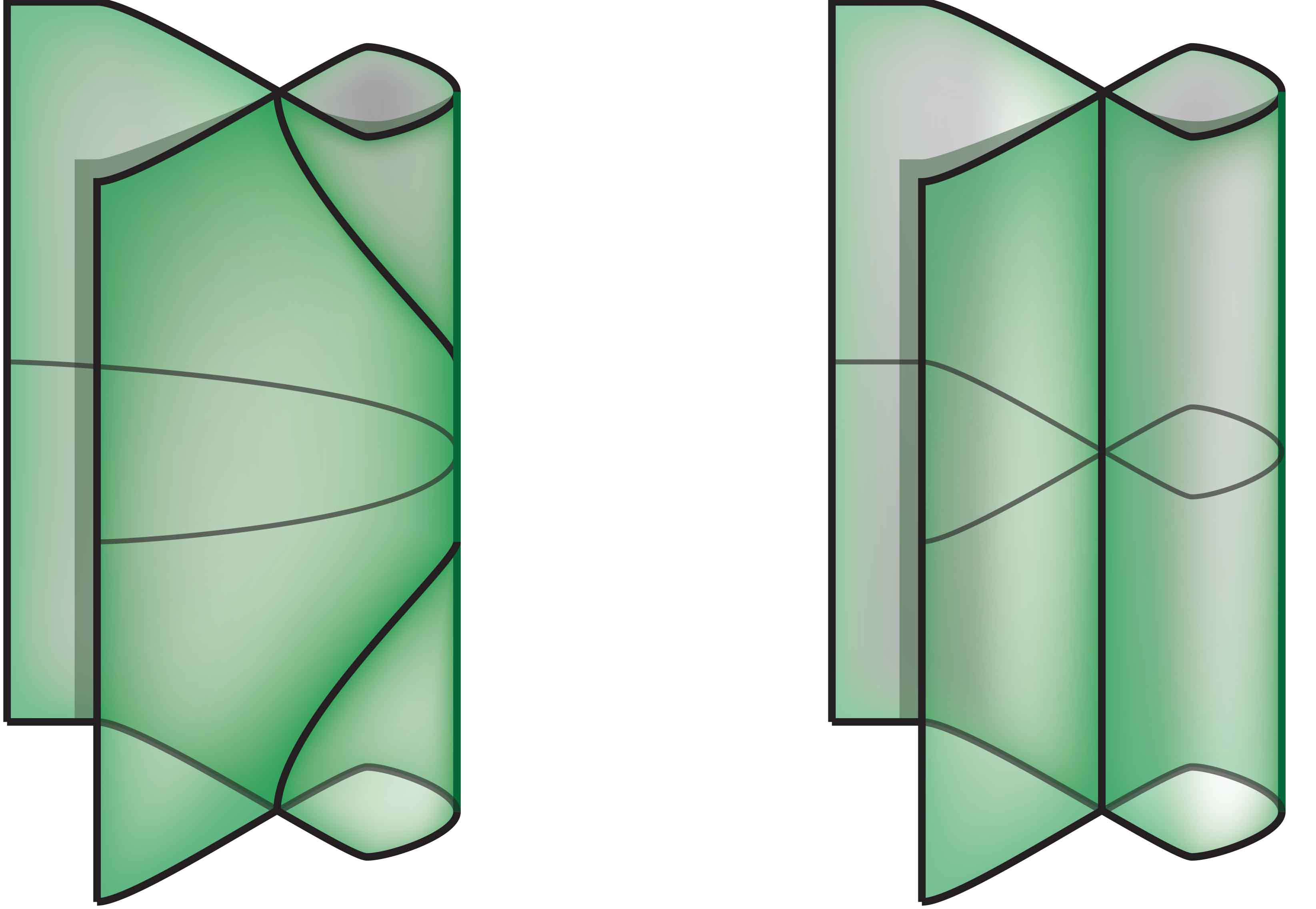}$ 
\caption{Fourth Roseman move $\mathcal{R}_4$: hyperbolic move of the type $\Omega_1$.} \label{ris:R4}
\end{center}
\end{figure}

\begin{figure}[h]
\begin{center}
$\includegraphics[width=90mm]{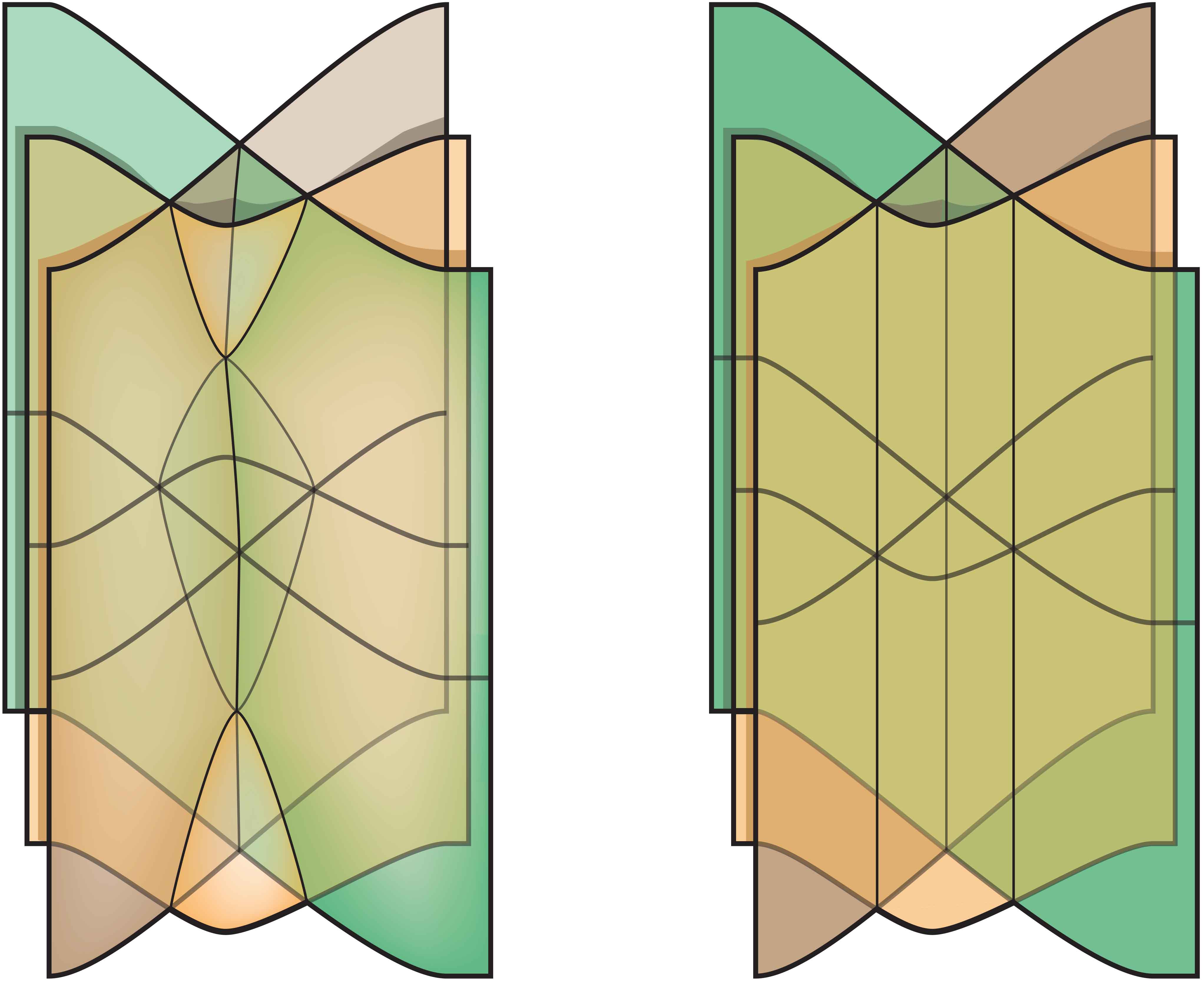}$ 
\caption{Fifth Roseman move $\mathcal{R}_5$: type $\Omega_3$.} \label{ris:R5}
\end{center}
\end{figure}

\begin{figure}[h]
\begin{center}
$\includegraphics[width=90mm]{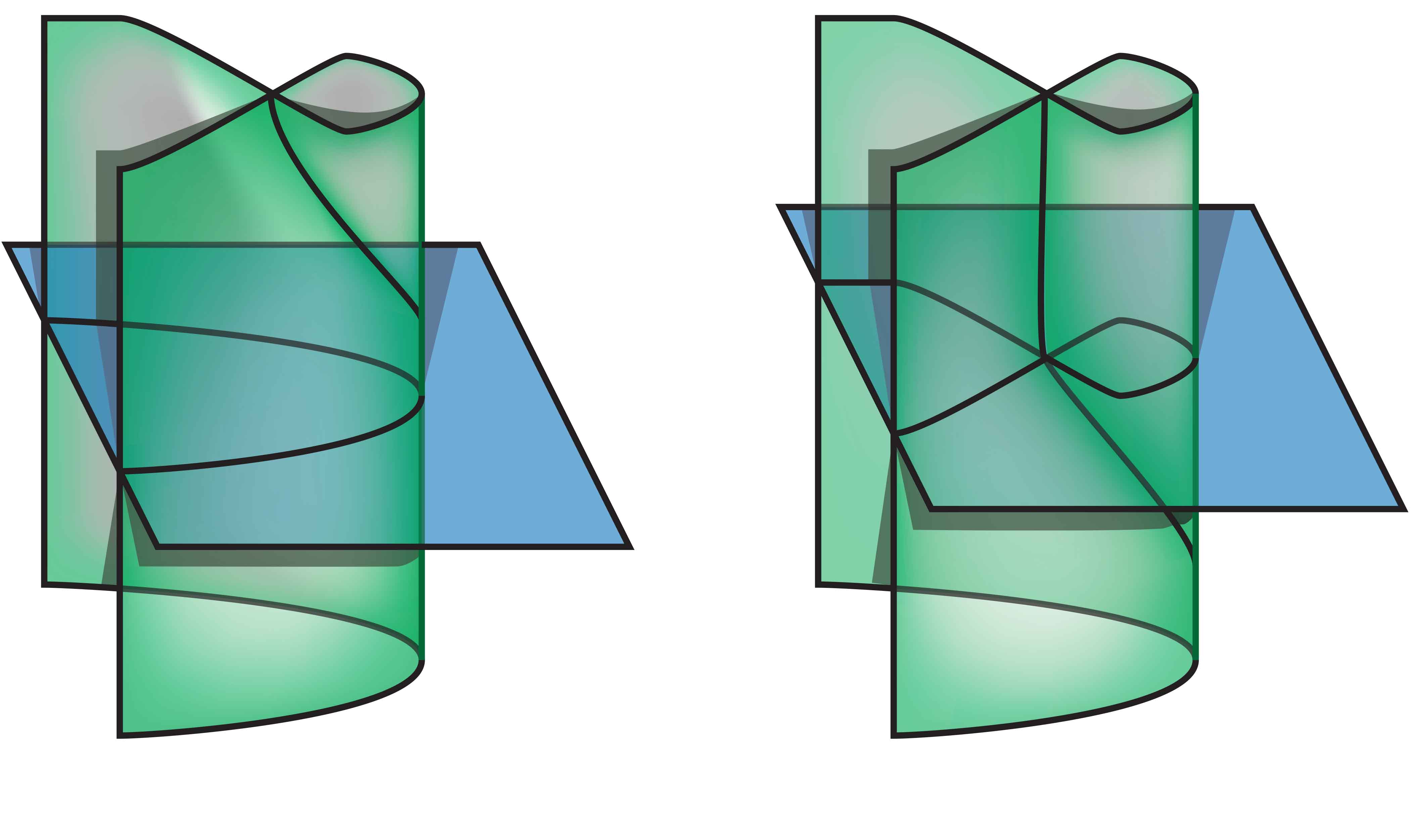}$ 
\caption{Sixth Roseman move $\mathcal{R}_6$: branch point crossing.} \label{ris:R6}
\end{center}
\end{figure}

\begin{figure}[h]
\begin{center}
$\includegraphics[width=90mm]{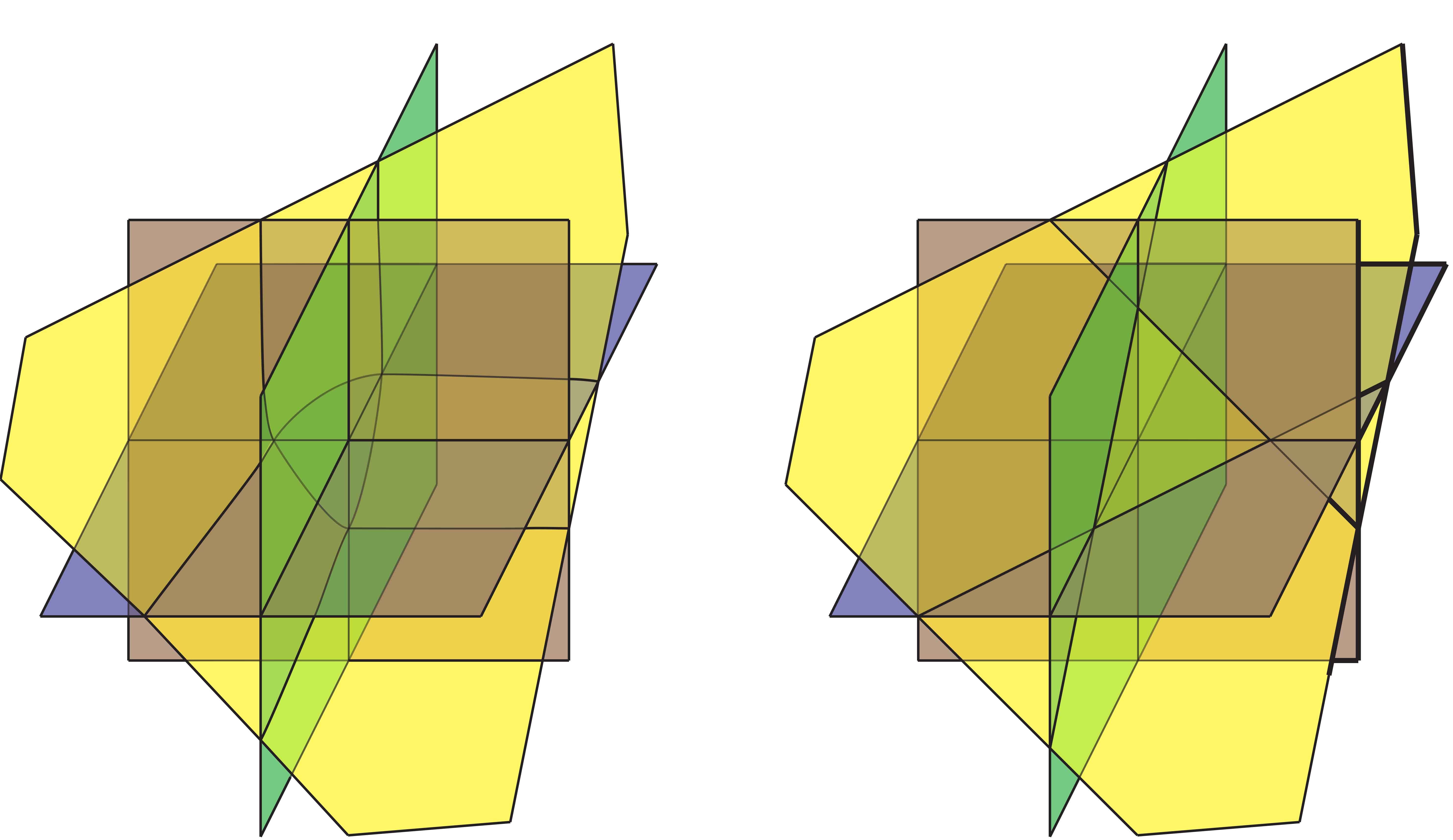}$ 
\caption{Seventh Roseman move $\mathcal{R}_7$: tetrahedral.} \label{ris:R7}
\end{center}
\end{figure}

Roseman \cite{ros} proved that two diagrams represented the same 2-knot (so their preimages differ by an isotopy) if and only if the diagrams can be connected by a finite sequence of isotopies of $\pi(L)\subset F$ and the so called {\it Roseman moves} moves shown in Fig. \ref{ris:R1}--\ref{ris:R7} (those pictures are originally due to Carter, see \cite{carter}). Figures \ref{ris:R1}--\ref{ris:R7} do not show which leaf is upper and which is lower. The corresponding moves take place for every ``over/under'' structure. In other words, the following statement holds:

\begin{teo}[Roseman \cite{ros}]
Two diagrams represent equal 2-knots if and only if they can be connected by a finite sequence of {\rm Roseman moves} shown in Figures \ref{ris:R1}--\ref{ris:R7} and isotopies of the diagram in $\mathbb{R}^3$.
\end{teo}

Consider a 2-knot diagram. It contains singularities of three types: double points, triple points and branch points (cusps). Triple points and cusps lie in the closure of the double point set. It is easy to see that this closure is an image of a one-dimensional manifold embedded into a 3-space. It is neat to include cusps into the set of double points. In that case the boundary of the double point manifold consists of cusps. Every triple point has exactly three preimages in the double points manifold.

Consider a 2-knot $K$. By definition it is an image of an embedding $p$ of a 2-sphere $S$ into $F\times [0,1]$. Consider its diagram $D$. It can be regarded as an image of the sphere $S$ under a map $\psi =\pi \circ p$. Consider the preimages of double points, triple points and cusps on the sphere under this map. The resulting complex will be called a {\it spherical diagram} of the 2-knot. Formally one could say that this complex is obtained from a sphere by factorisation by the following relation: two points are equivalent if they correspond to the same point on a double line. Spherical diagram is a natural generalization of the notion of a Gauss diagram of 1-knot and can be defined regardless of the 2-knot in 4-space. More precisely:

\begin{opr}
\label{def:gauss}
A {\rm sperical diagram} is a 2-complex consisting of a sphere $S$ and a set $D$ of marked curves on it such that:
\begin{enumerate}
\item Every curve is either closed or ends with a cusp, the number of cusps is finite;
\item Every curve of the set $D$ is paired with exactly one curve of that set, one of the paired curves is marked as {\rm upper}, both curves are oriented (marked with arrows) up to simultaneous orientation change;
\item Two curves ending in the same cusp are paired and both arrows either look towards the cusp or away from it;
\item If two curves intersect, the curves paired with them intersect as well (thus a triple point appears on the sphere $S$ three times).
\end{enumerate}
\end{opr}

Let us clarify the significance of the arrows marking the curves on a spherical diagram.

Like every two points (ends of a chord) of a Gauss diagram of a 1-knot correspond to a crossing, every pair of paired curves correspond to a double line of a 2-knot. But unlike the one-dimensional case, two curves can be glued in two ways differing by the curves' orientation. The arrows fix this uncertainty: the curves are glued with respect to the arrows. The third condition of the definition \ref{def:gauss} means that the two curves ending in the same cusp yield one double line and the cusp is glued with itself, not the opposite end of one of the curves.

It is clear from the definition that every 2-knot has a sperical diagram. On the other hand, not every spherical diagram corresponds to a ``real'' 2-knot just like not every chord diagram corresponds to a classical 1-knot. It is worth noting that due to the decoration of the curves on a spherical diagrams with arrows and ``over/under'' information there is a bijection between space and spherical diagrams of any 2-knot $K$.

Equivalent diagrams of classical 2-knots are connected by Roseman moves $\mathcal{R}_1$--$\mathcal{R}_7$. Those moves have sperical conterparts $\mathcal{S}_1$--$\mathcal{S}_7$. Let us describe them in more detail.

Let us call a {\it local link} a proper embedding in a general position of a disjoint union of $n$ discs into a ball $D$. Images of the discs will be called the {\it leaves} of the local link. Diagrams of local links are defined exactly as those of (global) links. {\it Local moves} --- moves on local links --- remove a set of discs and replace it with a new one. Roseman moves can be regarded as local. since every Roseman move doesn't change diagram away from a small sphere. Local links appearing in Roseman moves can have one leaf (moves $\mathcal{R}_3$, $\mathcal{R}_4$), two (moves $\mathcal{R}_1$, $\mathcal{R}_2$, $\mathcal{R}_6$), three (move $\mathcal{R}_5$) or four leaves (move $\mathcal{R}_7$).

Let us introduce the notion of {\it spherical diagram of a local $n-$component link}.

\begin{opr}
\label{loc_gauss}
A {\rm sperical diagram of a local $n-$component link} is a 2-complex consisting of a sphere $S$, a set of disjoint discs $U_i$, $i=1, \dots , n$ on it and a set $D$ of marked curves on it such that:
\begin{enumerate}
\item Every curve lies on one of the discs $U_i$;
\item Every curve is either closed or ends with a cusp or ends on the boundary of the corresponding disc, the number of cusps is finite;
\item Every curve of the set $D$ is paired with exactly one curve of that set, one of the paired curves is marked as {\rm upper}, both curves are oriented (marked with arrows) up to simultaneous orientation change;
\item Two curves ending in the same cusp are paired and both arrows either look towards the cusp or away from it;
\item If two curves intersect, then the curves paired with them intersect as well (thus a triple point appears on the sphere $S$ three times).
\end{enumerate}
\end{opr}

A move on a local link induces a move on its spherical diagram. On the other hand, consider a given spherical diagram and a move on it changing the diagram only inside a finite set of disjoint discs $U_1, \dots , U_n$ (though the set of curves $D$ and their topology can change). Such move will be called {\it local}. More formally: \\

 {\it local move consists of a removal of local link from a spherical diagram and glueing in another local link with the same boundary respecting the curves meeting the boundary}. \\
 
It is naturally identified with a move on a local link diagram. Every Roseman move (or its spherical analog) corresponds to two local links: its left-hand side and right-hand side.

\begin{figure}[h]
\begin{center}
$\includegraphics[width=90mm]{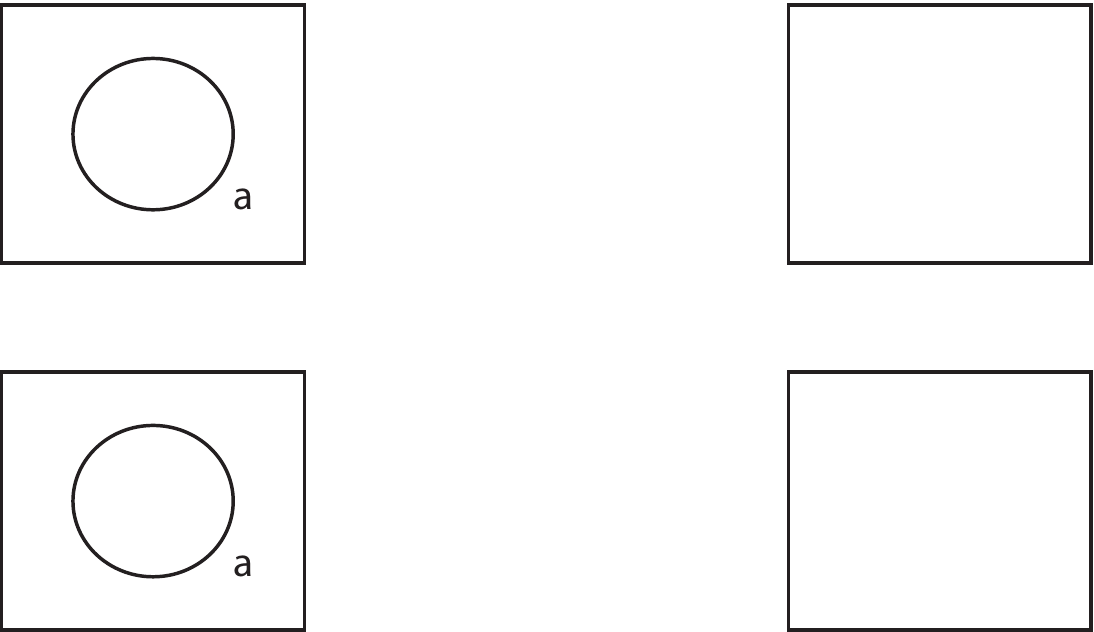}$ 
\caption{First spherical move $\mathcal{S}_1$.} \label{ris:S1}
\end{center}
\end{figure}

\begin{figure}[h]
\begin{center}
$\includegraphics[width=90mm]{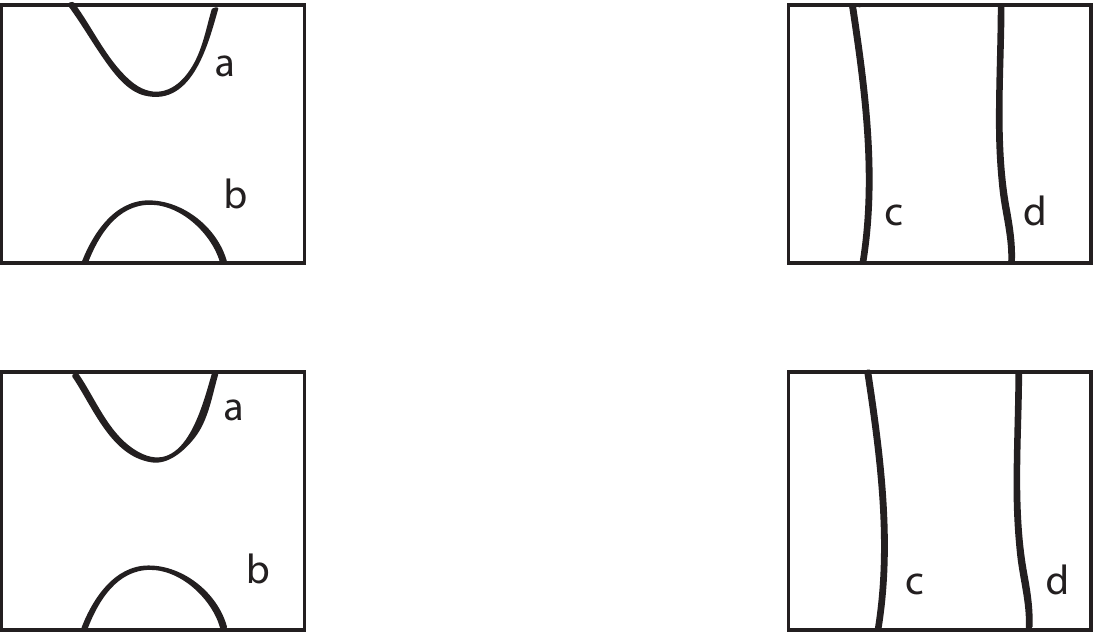}$ 
\caption{Second spherical move $\mathcal{S}_2$.} \label{ris:S2}
\end{center}
\end{figure}

\begin{figure}[h]
\begin{center}
$\includegraphics[width=90mm]{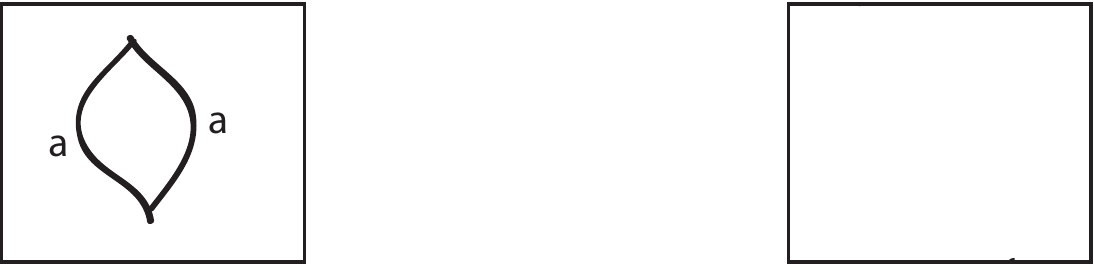}$ 
\caption{Third spherical move $\mathcal{S}_3$.} \label{ris:S3}
\end{center}
\end{figure}

\begin{figure}[h]
\begin{center}
$\includegraphics[width=90mm]{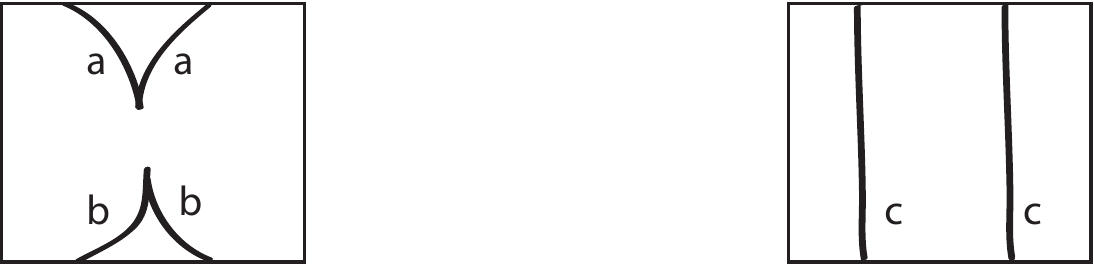}$ 
\caption{Fourth spherical move $\mathcal{S}_4$.} \label{ris:S4}
\end{center}
\end{figure}

\begin{figure}[h]
\begin{center}
$\includegraphics[width=90mm]{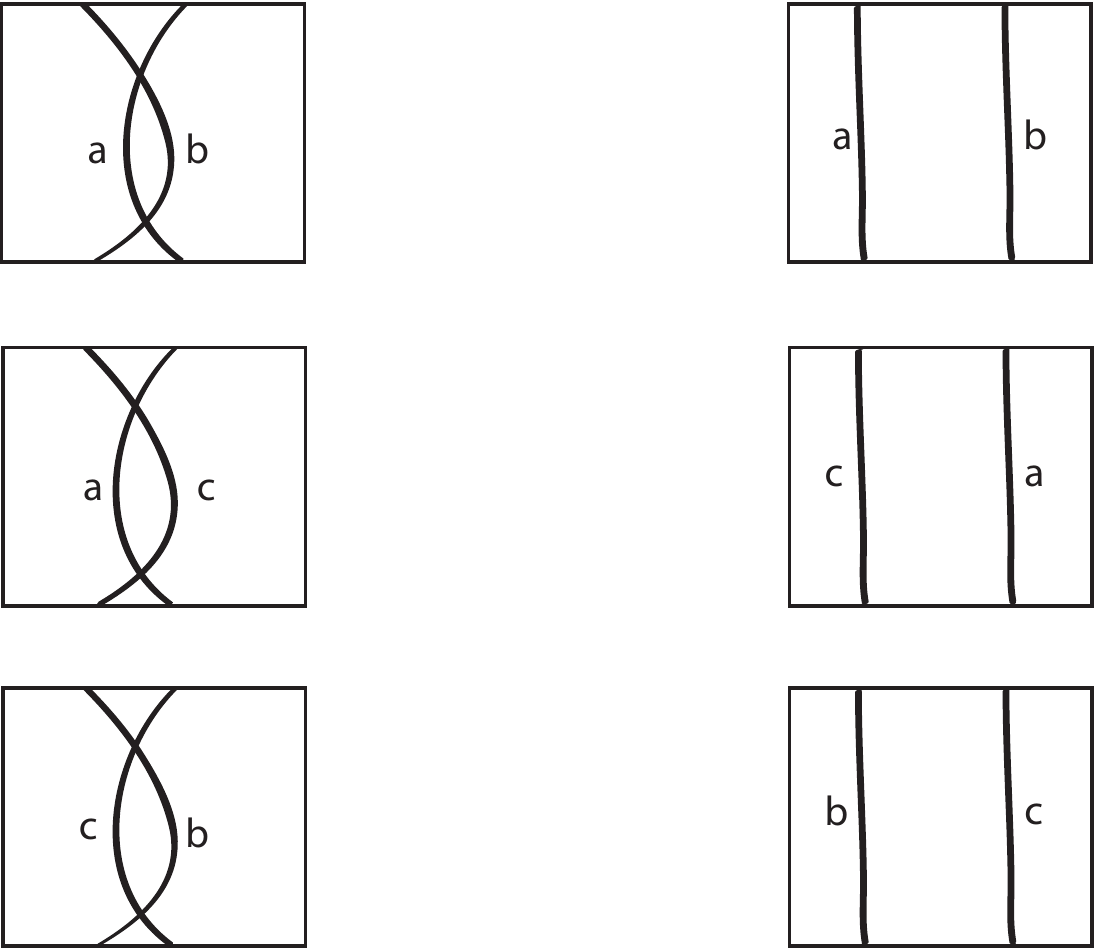}$ 
\caption{Fifth spherical move $\mathcal{S}_5$.} \label{ris:S5}
\end{center}
\end{figure}

\begin{figure}[h]
\begin{center}
$\includegraphics[width=90mm]{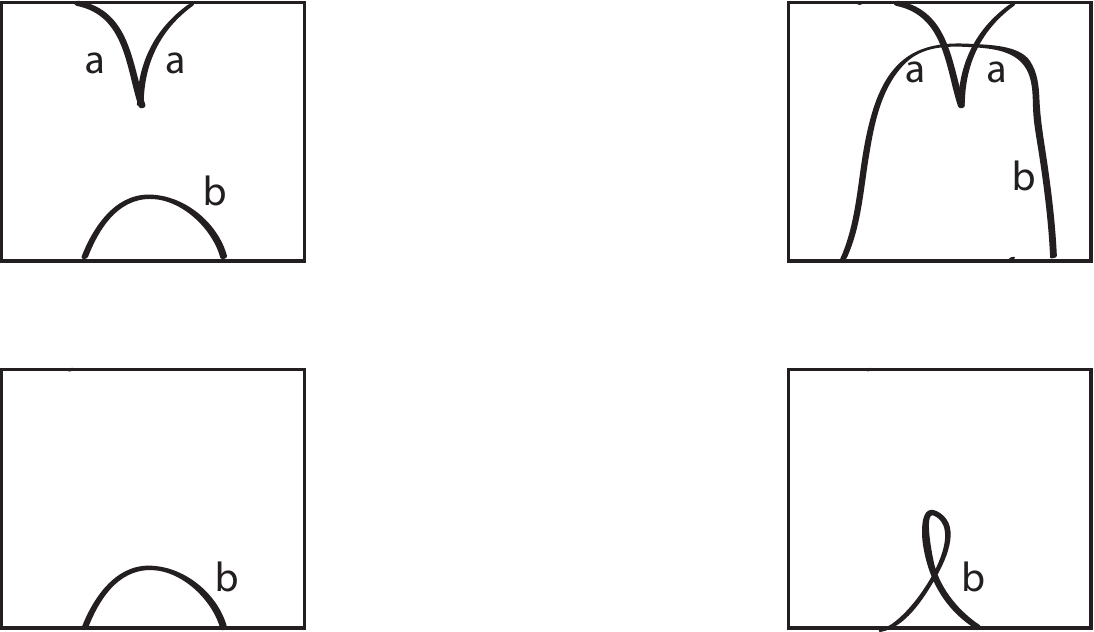}$ 
\caption{Sixth spherical move $\mathcal{S}_6$.} \label{ris:S6}
\end{center}
\end{figure}

\begin{figure}[h]
\begin{center}
$\includegraphics[width=90mm]{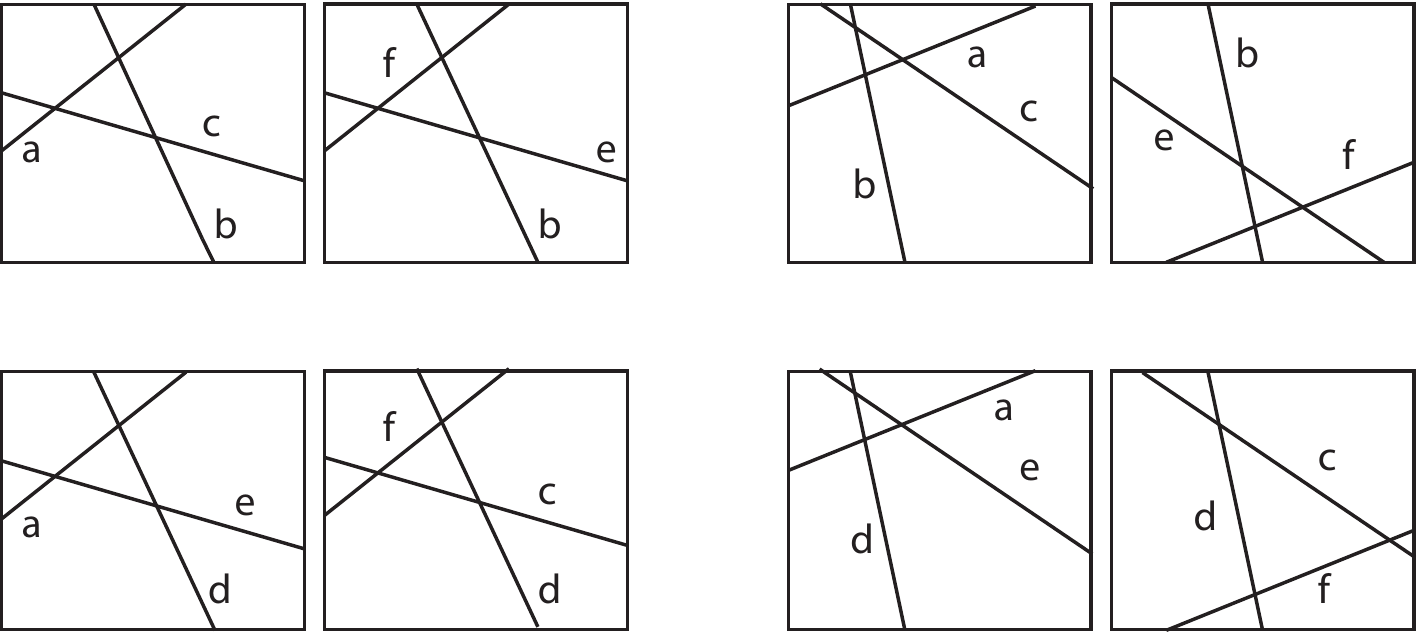}$ 
\caption{Seventh spherical move $\mathcal{S}_7$.} \label{ris:S7}
\end{center}
\end{figure}

Thus, since every Roseman move $\mathcal{R}_i$ is a local move, it corresponds to a (local) move on a spherical diagram $\mathcal{S}_i$. Those moves are depicted in Figures \ref{ris:S1} -- \ref{ris:S7}. The same letters denote paired curves; the same letters in the left-hand and right-hand sides of the moves denote the curves in a natural one-to-one correspondence.

Note that the definition of the moves $\mathcal{S}_i$ does not require the spherical diagram to be a diagram of a ``real'' 2-knot. We introduce the following definitions:

\begin{opr}
An {\rm abstract 2-knot} is an equivalence class of spherical diagrams modulo moves $\mathcal{S}_1 - \mathcal{S}_7$.
\end{opr}

\begin{opr}
A {\rm free 2-knot} is an abstract 2-knot without the ``over/under'' structure in double and triple points.
\end{opr}

Additionally, note the following fact. Every Roseman move is done as follows: one cuts a local link out and glues a different local link in. Therefore, if a spherical diagram was a diagram of a real 2-knot, the resulting diagram also corresponds to a real 2-knot.

\section{Parities on 2-knots and links}

The notion of {\it parity} was originally devised by V.O. Manturov (see for example \cite{cob}) for 1-knots. Its gist is to devide crossing into {\it even} and {\it odd} so that some natural conditions on the parity behavior under the Reidemeister moves hold.

In case of 1-knots the simpliest --- and in some sense the most important --- parity (Gaussian) is constructed from a chord diagram of a knot. Recall, that in case of the Gaussian parity a chord (and the corresponding crossing) was defined as even if it intersected an even number of chords. For classical knots the parity is always trivial (all crossings are even); for virtual knots it is not so.

Apart from the Gaussian parity there are many parities in virtual knot theory. It is possible to define a set of {\it parity axioms} so that for every parity satisfying them --- regardless of its individual properties --- one can construct invariants and prove theorems obtained before for a given parity only (say, the Gaussian).

In the present paper the same is done for 2-knots. In one-dimensional case crossings were named even and odd. In two-dimensional case a parity is defined for the set of {\it double lines}. First of all we define the Gaussian parity (as an analog to the 1-dimensional Gaussian parity) and then the general principle of parity will be formulated.

\subsection{The notion of Gaussian parity}

First of all let us define the Gaussian parity following \cite{cob}. We begin with the definition of the Gaussian parity for a real 2-knot and then generalize it to abstract 2-knots.

Consider a 2-knot $K$ and its diagram $L$ which is an embedding of a sphere into $\mathbb{R}^3$. Now consider its spherical diagram $S$.

Consider a double point $x\in D_0$. It has two preimages and, therefore, two corresponding points on the diagram $S$; denote them by $x_1, x_2$.

Consider an arbitrary loop $\eta$ connecting $x$ to $x$ whose preimage $\tilde{\eta}$ is in general position with respect to $D^*$: we suppose all intersections between $\tilde{\eta}$ and $D^*$ to be transverse. All such curves $\tilde{\eta}$ lie on the sphere $S$ and have fixed endpoints, thus homotopic. We only need to fix the behaviour of the curve near $x_1$ and $x_2$.

Consider a double line $\zeta$ which contains the point $x$. Orient it arbitrarily and orient its preimages $\zeta_1 \cap U(x_1)$ and $\zeta_2 \cap U(x_2)$ in the coherent way. We impose the following condition: bases $(\dot{\zeta_1}, v_1)$ and $(\dot{\zeta_2}, v_2)$ give different orientation to the sphere $S$; $\dot{\zeta_i}$ denotes the unit tangent vector to the curve $\zeta_i$.

Now we call the parity of the point $x$ the parity of the number of elements of the set $Card(\tilde{\eta}\cap D)$. Informally one can say that we count the number of double lines intersecting the curve connecting two preimages of a double point --- provided the ends of the curve are properly oriented.

Note that simultaneous change of orientation of the preimages of the double line induced by the orientation change of the double line $\zeta$ itself doesn't change the parity defined this way.

Thus we defined parity of an arbitrary double point $x$. We call it {\it the Gaussian parity of the double point}.

Elementary check verifies the following statement (see \cite{cob}):

\begin{state}
The Gaussian parity is constant along a double line.
\end{state}

Now we can define the Gaussian parity of a double line:

\begin{opr}
Consider an arbitrary double point $x\in D_0$ on a double line $\gamma$ and consider its two preimages $x_1, x_2$. Connect the point $x_1$ with the point $x_2$ by a path $\tilde{\eta}$ so that the behaviour of $\gamma$ near the endpoints is compatible (as before). {\rm the Gaussian parity} of the double line $\gamma$ is the parity of the number of intersections between $\tilde{\eta}$ and the preimage of the set $D^*$.
\end{opr}

Now we can define the Gaussian parity for abstract and free 2-knots. Note that the Gaussian parity of a double point (and thus that of a double line) was defined solely by using a spherical diagram. Since abstract (and free) knots are equivalence classes of spherical diagrams (maybe with additional information encoded) modulo moves, the definition of the Gaussian parity for them can be given verbatim. \\

Due to above mentioned considerations there is sense in the study of the parity of the double lines appearing in Roseman moves: it is interesting to understand what can be said of thier parity regardless of the structure of the rest of a knot. We do not specify in which space the knot lies, it is considered abstractly (as a local link, to be precise). Despite the fact that the Gaussian parity is trivial for classical knots, the parity can be non-trivial for abstract knots. And since every transformation of a spherical diagram is (locally) a Roseman move, it is important to study the parity of the double lines involved. \\

First of all note that the parity of double lines not involved in a move remains unchanged.

Now consider the Roseman moves 1--7 and the corresponding local links.

\begin{enumerate}
\item Moves 4 and 6 yield that every double line ending in a cusp is even.
\item Two double lines from the second move (located ``closely'') have the same the Gaussian parity.
\item The double line in the right-hand side of the fourth move is even.
\item The third move creates a double line; it is always even.
\item Among the lines meeting at a triple point there is an even number of odd ones.
\item There are three double lines in the fifth move. There are either two o zero odd among them.
\end{enumerate}
 
Note that the fifth, sixth and seventh moves do not change the number of double lines. The double lines before and after the move are in a natural one-to-one correspondence. Parities of the corresponding lines are the same. \\

Later we will see that those properties are not random and stem from general axioms of 2-parity which we deduce from general approach described in the next subsection and are closely related to parities on 2-links.

\subsection{General parity principle}

Consider a more general situation: a disjoint union of $n-$spheres $S^n$ in $\mathbb{R}^{n+2}$. Consider them from a {\it diagrammatic} viewpoint: as a projection in general position to a subspace $\mathbb{R}^{n+1}$. Roseman proved that in any dimension there is a {\it finite} set of moves with the following property: two diagrams represent isotopic links if and only if they can be connected by a finite sequence of diagrams such that each next diagram is obtained using one of those moves or a trivial isotopy (see, for example, \cite{ros2}). For low dimensions those moves are completely described (that is, Reidemeister moves and Roseman moves).

Now we define {\it parity axioms} corresponding to the moves.

Consider a move $M$ on a diagram $D$ of a local link sending it to a diagram $D'$. There is a natural bijection between the leaves of the diagram $D$ and those of the diagram $D'$. Colour the leaves of the diagram $D$ arbitrarily with two colours. colour the corresponding leaves of the diagram $D'$ with the same colours (the compatibility condition can limit the colourings choice). Now we call a codimension 1 set {\it even} it is the intersection set of two leaves of the same colour and {\it odd} otherwise. Finally to define a {\it parity axiom corresponding to the move $M$} we take a condition on the parity of the intersection sets which holds for all possible colourings of this local link. Note that different moves can possibly yield the same axiom. \\

\begin{figure}[h] 
\begin{center}
$\includegraphics[width=90mm]{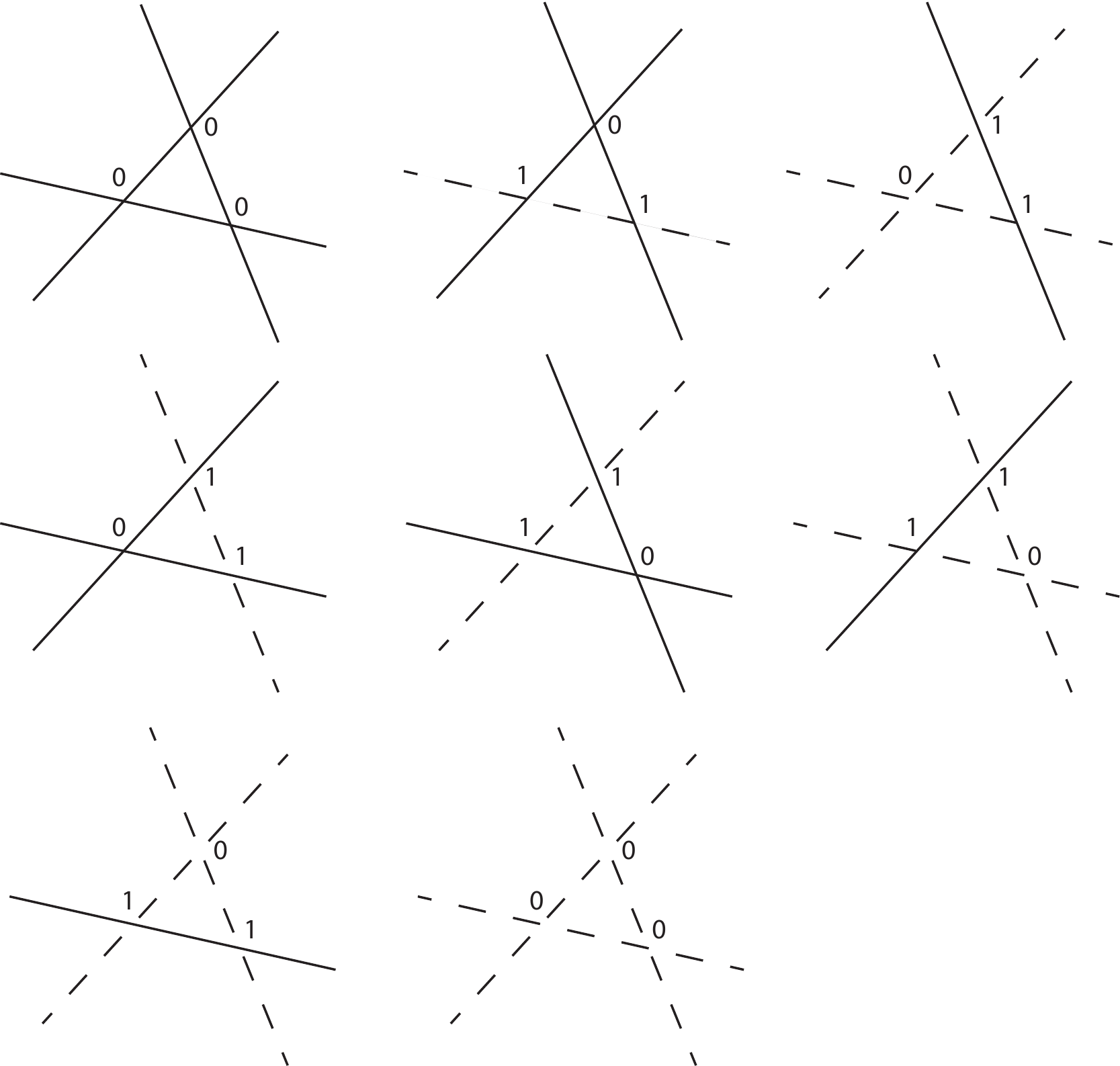}$ 
\caption{Possible colourings and parity of crossings of the third Reidemeister move.} \label{ris:reid3colour}
\end{center}
\end{figure}

\textbf{A clarifying example: the third Reidemeister move.} Local link which is changed by this move consists of three arcs (we don't care about over- and undercrossings for our purposes). Each can be coloured in one of two colours. That gives eight possible colourings of the local link (see Fig. \ref{ris:reid3colour}). Every colouring yields a parity of each of the three crossings. But there is a property which holds for all eight colourings: the sum of parities modulo 2 equals zero. That is the parity axiom corresponding to the move $\Omega_3$. \\

In the case of 2-knots and 2-links the general principle can be used as is since all the moves are known: those are Roseman moves. They yield the following eight parity axioms for double lines:

\begin{enumerate}
\item[0.] Parity is constant along double lines;
\item Double line ending in a cusp is even;
\item The sum of parities of double lines meeting in a triple point equals $0 \mod 2$;
\item All double lines in the move $\mathcal{R}_2$ have the same parity;
\item Among the three double lines from the move $\mathcal{R}_5$ there are either two or zero odd ones;
\item The double line appearing in the move $\mathcal{R}_3$ is even;
\item There is a natural bijection between the lines on the left- and on the right-hand sides of the moves $\mathcal{R}_5, \mathcal{R}_6, \mathcal{R}_7$ induced by the correspondence of the leaves of the diagrams. Parity is preserved by the correspondence;
\item The double line from the right side of the move $\mathcal{R}_4$ is even.
\end{enumerate}

Besides, it is always supposed that the moves do not change the parity of any double lines not present in the local link. \\

Note though, that those axioms --- obtained from the moves --- are redundant. Indeed, the axiom 5 follows from the axiom 1 since the line in question ends in a cusp. The axiom 4 follows from the axioms 2 and 6 since in the left-hand side of the move there are two triple points.

The axioms 3 and 7 seem independent of others. But they can be avoided by using the following trick.

Every double line of a local link is either closed or meets the boundary twice. We will call the connected components of the intersection set between the boundary and the double lines {\it boundary double points}. There is a natural bijection between the leaves of the local links in the left and the right sides of Roseman moves. It induces the bijection between their boundaries which respects the boundary double points. Therefore, though there is not always a bijection between the double lines, one can always find a bijection between boundary double points. 

Thus we can reformulate the axiom 6 in the following way: \\

{\it (correspondence axiom) the corresponding boundary double points have the same parity.} \\

This axiom together with axioms 0 and 1 yield axioms 3 and 7. Indeed, the left side of the fourth move have two double lines ending with cusps. Due to axiom 1 they are even (and the corresponding boundary double lines are even as well). Therefore due to the correspondence axiom the double line in the right side is even as well. In the same way one can check that the double lines in the second move have the same parity.

In that manner we have a smaller set of parity axioms. Note that they are less dependent of the actual list of Roseman moves: only the correspondence axiom appeals to them.

\begin{enumerate}
\item {\it Continuity axiom.} Parity is constant along double lines.
\item {\it Selfcrossing axiom.} A double line ending with a cusp is even.
\item {\it Triple point axiom.} The sum of parities of three double lines meeting in a triple point equals $0 \mod 2$.
\item {\it Correspondence axiom.} There is a natural bijection between boundary double point on the left and right sides of a Roseman move induced by the bijection of the diagram leaves. The parities of the corresponding points are the same.	
\end{enumerate}

\begin{opr}
Let $\mathcal{L}$ be a class of 2-lniks in $\mathbb{R}^3$, and let $A$ be the set of double line of their diagrams. A mapping $P: A \to \mathbb{Z}_2$ is called {\rm parity} if it satisfies the axioms 1 -- 4.
\end{opr}

As one can see, the axioms 0 -- 7 above are exactly the same as the properties we got when studying the Gaussian parity. Since the two systems of axioms are equivalent, we got the following

\begin{state}
The Gaussian parity is a parity.
\end{state}

Using the similar techniques one can define parity axioms for links in arbitrary dimensions.

\subsection{The notion of weak parity}

Besides parity there is a similar notion of {\it weak parity}. In the one-dimensional case it was developed, say, in \cite{nik}. We will describe it in an arbitrary dimension.

Since the number of moves in any given dimension is finite, denote by $N$ the maximal number of leaves of a local link from those moves. Now consider a move $M$ on a local link diagram $D$ sending it to a diagram $D'$. 

Colour the leaves of the diagram $D$ arbitrarily with $N$ colours. colour the corresponding leaves of the diagram $D'$ in the same colours (the compatibility condition can limit the colourings choice). Now we call a codimension 1 set {\it even} it is the intersection set of two leaves of the same colour and {\it odd} otherwise. 

Finally to define a {\it weak parity axiom corresponding to the move $M$} we take a condition on the parity of the intersection sets which holds for all possible colourings of this local link. Note that different moves can possibly yield the same axiom.

As we can see, the definition of weak parity is essentially the same as before, but the set of colours is larger (allowing all leaves in any move to have different colours). \\

As an example consider once again the third Reidemeister move. It is easy to verify that the three-colour colourings give only one new case: when all the leaves have different colour. Thus we get an axiom of (one-dimensional) weak parity: among the crossing in the third Reidemeister move there are 0, 2 or 3 odd ones, see \cite{nik}. \\

In two-dimensional case one has Roseman moves. The maximal number of leaves is 4 and is attained on the seventh move. The moves $\mathcal{R}_1, \mathcal{R}_2, \mathcal{R}_3, \mathcal{R}_4, \mathcal{R}_6$ deal with two or less leaves so the weak parity axioms for them are the same as (strong)  parity axioms. Now we study the axioms coming from the moves $\mathcal{R}_5$ and $\mathcal{R}_7$.

The move $\mathcal{R}_5$ includes 3 leaves so the only new colouring is the colouring of all leaves in different colours. In that case all double lines are odd. Therefore the corresponding weak parity axiom is: \\

 {\it there is a bijection between the double lines of the two local links appearing in the fifth Roseman move; among the lines there are 0, 2 or 3 odd ones and the corresponding lines have the same parity}.

The move $\mathcal{R}_7$ includes four leaves so there are two possible new colourings: {\it a)} all leaves have different colours, {\it b} three leaves are coloured in different colours and the fourth colour coincides with one of them. There are 6 double lines in the move, so in case a) their sum is still 0 mod 2. It is easy to check that in case b) one has 5 odd lines. Thus the corresponding weak parity axiom is: \\

{\it there is a bijection between the double lines of the two local links appearing in the seventh Roseman move; among the lines there are 0, 2, 4, 5 or 6 odd ones and the corresponding lines have the same parity}. 

Note that in case of weak parity there is no triple point axiom. \\

From the definition we get

\begin{state}
Every parity is a weak parity. 	
\end{state}

It is worth noting that both parity and weak parity axioms were defined from local links and colouring if their components. In the next section we discuss connections of this approach to actual links.

\subsection{Parities for 2-links}

Consider a 2-component 2-link. A {\it componentwise parity} can be defined in that case as for 1-links. There are double lines of two types: selfcrossing lines and mixed crossing lines. We call the the former lines even and the latter --- odd. 

\begin{state}
That decoration of double lines is a parity.
\end{state}

{\bf Proof.} We need to check the axioms 1 -- 4. It is easier to check the equivalent axioms 0 -- 7 since they come directly from moves. Indeed, the link components are naturally ``coloured'' into two colours (one for each component). Thus by construction the axioms are satisfied. The proof is concluded. \\

This is not a coincidence. When defining the parity axioms we considered a local link. But in that case those local links can be continued to the actual link we consider. Thus the axioms are satisfied automatically.

From that point of view the componentwise parity is the most basic ones since the axioms stem form links. \\

Now consider a multicomponent link.

First of all we can call the selfcrossing double lines even and all other --- odd. For the same reasons as before we get a {\it weak} parity (again the most basic one). It is called componentwise (weak) parity as well.

Moreover a {\it multivalued parity} or a {\it group-valued parity} can be defined. A multivalued parity is such a mapping from the set of double lines of a 2-knot or link into a finite commutative group $G$ that the sum of parities of all double lines in a local link from a Roseman move is zero in the group $G$. From that point of view the (strong) parity defined above is a parity valued in $\mathbb{Z}_2$.

\subsection{Parity projection and hierarchy}

At the level of free 2-knots (and therefore spherical diagrams) one can define a {\it parity projection} like in one-dimensional case (see for example \cite{knot_parity}). It is done as follows.

Denote by $\mathcal{G}$ the set of spherical diagrams of a given knot. Define a mapping  $p: \mathcal{G} \to \mathcal{G}$ in the following way: the image of a spherical diagram under the projection is a spherical diagram, obtained by deleting all odd curves from the set $D$. The following theorem holds:

\begin{teo}
\label{theo_proj}
The mapping $p$ is well defined. That is, consider free 2-knots $K$ and $K'$ and their spherical diagrams $D(K)$ and $D(K')$. If the knots $K$ and $K'$ are equal, then the knots $\tilde{K}$ and $\tilde{K'}$ are equivalent as well, where the knots are such that $D(\tilde{K})=p(D(K)), \, D(\tilde{K'})=p(D(K'))$.
\end{teo}

The proof of this theorem directly follows from Roseman moves (or their spherical analogs).

The projection mapping leads to the notion of {\it parity hierarchy}. Consider a spherical diagram of a knot $K$ with parity. Project it by using the mapping $p$. We get a new diagram of some new 2-knot. Its double lines again can be named even and odd according to the same rules as before. Note that the lines which used to be even can now be odd.

This process can be continued. For any abstract 2-knot $K$ there exists a natural number $n$ such that $p^{n+1}(D(K)) = p^n(D(K))$. In other words, the parity hierarchy stabilizes. We will use the notation $p^{stab}=p^n$. Theorem \ref{theo_proj} yields

\begin{teo}
\label{theo_stab}
The mapping $p^{stab}$ is well defined.
\end{teo}

The notion of projection is repeated verbatim in the weak parity case. We will denote the corresponding mappings $p_w$ and $p_w^{stab}$. As before we have the following theorems:

\begin{teo}
\label{weak_proj}
The mapping $p_w$ is well defined. That is, consider free 2-knots $K$ and $K'$ and their spherical diagrams $D(K)$ and $D(K')$. If the knots $K$ and $K'$ are equivalent, then the knots $\tilde{K}$ and $\tilde{K'}$ are equivalent as well, where the knots are such that $D(\tilde{K})=p_w(D(K)), \, D(\tilde{K'})=p_w(D(K'))$.
\end{teo} 

\begin{teo}
\label{weak_stab}
The mapping $p_w^{stab}$ is well defined.
\end{teo}

\subsection{The succession principle}

The notion of parity projection is closely related to the following important principle: the {\it succession principle}. Consider a set of objects with moves on them; the objects are considered equal if they can be connected by a sequence of those moves. Then consider a decoration of those objects: marking the co-dimension 1 starts with some additional information and the sufficient adaptation of the moves.

A classical example of such theory is knot theory. The initial objects are knot diagrams which can be decorated in various ways --- making some crossings virtual, introducing parity, etc. The moves in that case are Reidemeister moves and their virtual and other analogs.

The set of objects without decoration can be naturally included in the set of decorated objects. Indeed, every object can be regarded as trvially decorated. The {\it succession principle} is the following: \\

{\it if two objects regarded as decorated with trivial markings are equivalent as decorated objects, they are equivalent in the class of objects with trivial decoration}. \\

An example of this principle in action is given by a well known theorem: two classical knots equivalent as virtual are equivalent as classical.

In every given case the proof of this principle could be obtained by some techniques similar to parity projection method.

To illustrate it we prove the following theorem:

\begin{teo}
Consider an arbitrary weak parity $p$. Consider two free 2-knots $K_1, K_2$ for which the parity $p$ is trivial. Let those knots be equal in the class of knolts with the parity $p$. Then they are equal in the class of knots for which the parity $p$ is trivial.
\end{teo}

{\bf Proof.} Consider the diagrams $D_1=D(K_1)$ and $D_2=D(K_2)$. Since the knots $K_1$ and $K_2$ are equal, there exists a finite sequence of diagrams $D_1=D^0\sim D^1 \sim \dots \sim D^n = D_2$, and for every $i$ the diagrams $D^i$ and $D^{i+1}$ are equivalent via a Roseman move or a trivial isotopy.

Due to Theorem \ref{weak_stab} we have $\tilde{D}^0\sim \dots \sim \tilde{D}^n$, where $\tilde{D}^i=p_w^{stab}(D^i)$, and $\tilde{D}^i$ and $\tilde{D}^{i+1}$ are again equivalent via a Roseman move or a trivial isotopy.

On the other hand, since the weak parity is trivial on the diagrams $D_1$ and $D_2$ we have $p_w^{stab}(D_1)=D_1$ and $p_w^{stab}(D_2)=D_2$. Thus the diagrams are connected by a sequence of diagrams with weak parity $p$ obtained from one another via Roseman moves and trivial isotopy.

\section{Parity in 2-knot invariants}

The notion of parity allows one to refine existing invariants of knots and links and to create new ones. Here are several examples. \\

{\bf 1. The number of odd triple points.} The first example will be the following.

Earlier in this paper we attributed parity to codimension 1 singularities --- double lines. But triple points can also be decorated with parity in a natural way.

Fix a weak parity $P$ and consider a triple point $x$. Three double lines  $\alpha, \beta$ and $\gamma$ meet at that point and have weak parities $P(\alpha), P(\beta), P(\gamma) \in \{0,1\}$. We call the {\it triple point $x$ parity} the number $P(x):=P(\alpha)+P(\beta)+P(\gamma) \mod 2$.

An elementary check proves the following theorem:

\begin{teo}
The parity of the number of triple points of a diagram is an invariant of (free) 2-knots.
\end{teo}

{\bf Example.} To illustrate the use of the described invariant, let us consider the following class of 2-knots: {\it spun knots}. A spun knot is constructed from a free 1-knot. Consider a free 1-knot $K$ and take its Gauss diagram. Take a diameter of the diagram a spin it around the chosen diameter. The diagram will thus yield a sphere with parallels spun by the ends of the chords. It is always possible to choose the diameter in such a way that different endpoints yield different parallels. Those parallels are paired in the same way as the endpoints the came from. Thus one obtains a spherical diagram of a free 2-knot. That knot is a spun knot $S(K)$.

It is easy to verify, that the construction of a spun knot is well defined, that is, equivalent 1-knots yield equivalent 2-knots.

Parity on spun knots is inherited from the underlying 1-knots.

Now we can give a simple example of the use of the invariant. Consider an arbitrary free 2-knot $K$ with an odd number of (weakly) odd triple points. Since any spun knot has no triple points, the knot $K$ is not spun; in other words, there is no 1-knot $K'$ such that $K = S(K')$. \\

{\bf 2. 2-quandles.} Another knot invariant which can be refined via parity is {\it quandle} (see, say, \cite{winter}). In 1-dimensional case the refinement was done by D.M. Afanasiev. Let us describe this invariant in more detail.

First of all recall the definition of a quandle due to \cite{winter} (see also \cite{Joyce, Fenn, Matveev}).

\begin{opr}
A {\rm quandle} is a set $Q$ with an operation $\circ$, satisfying the following conditions:
\begin{enumerate}
    \item $a\circ a = a$ for all $a\in Q$;
    \item For all $a,b\in Q$ there exists a unique $x\in Q$ such that $x\circ b=a$;
    \item $(a\circ b)\circ c = (a\circ c)\circ (b\circ c)$ for all $a,b,c \in Q$.
\end{enumerate}
\end{opr}

The second condition (the existence of an inverse) lets us to say that there is an inverse operation denoted with $\bar{\circ}$. 

The simpliest example of quandle operatiion is given by a conjugation in case of a set $Q$ with group structure. In other words, one can put $a\circ b = bab^{-1}$. It is easy to see that all the quandle axioms are satisfied in that case.

Quandles are closely related to knots and their invariants. An abstract knot diagram gives rise to a quandle in the following way. Consider a set $Q$ elements of which are in a one-to-one correspondence with connected components of the diagram (here we consider the leaves to be broken at codimension 1 singularities according to the ``over/under'' structure). At each codimension 1 strat there are three connected components. We use conjugation as the quandle operation. Thus we get a {\it free quandle}  $(Q,\circ)$.

\begin{figure}[h] 
\begin{center}
$\includegraphics[width=90mm]{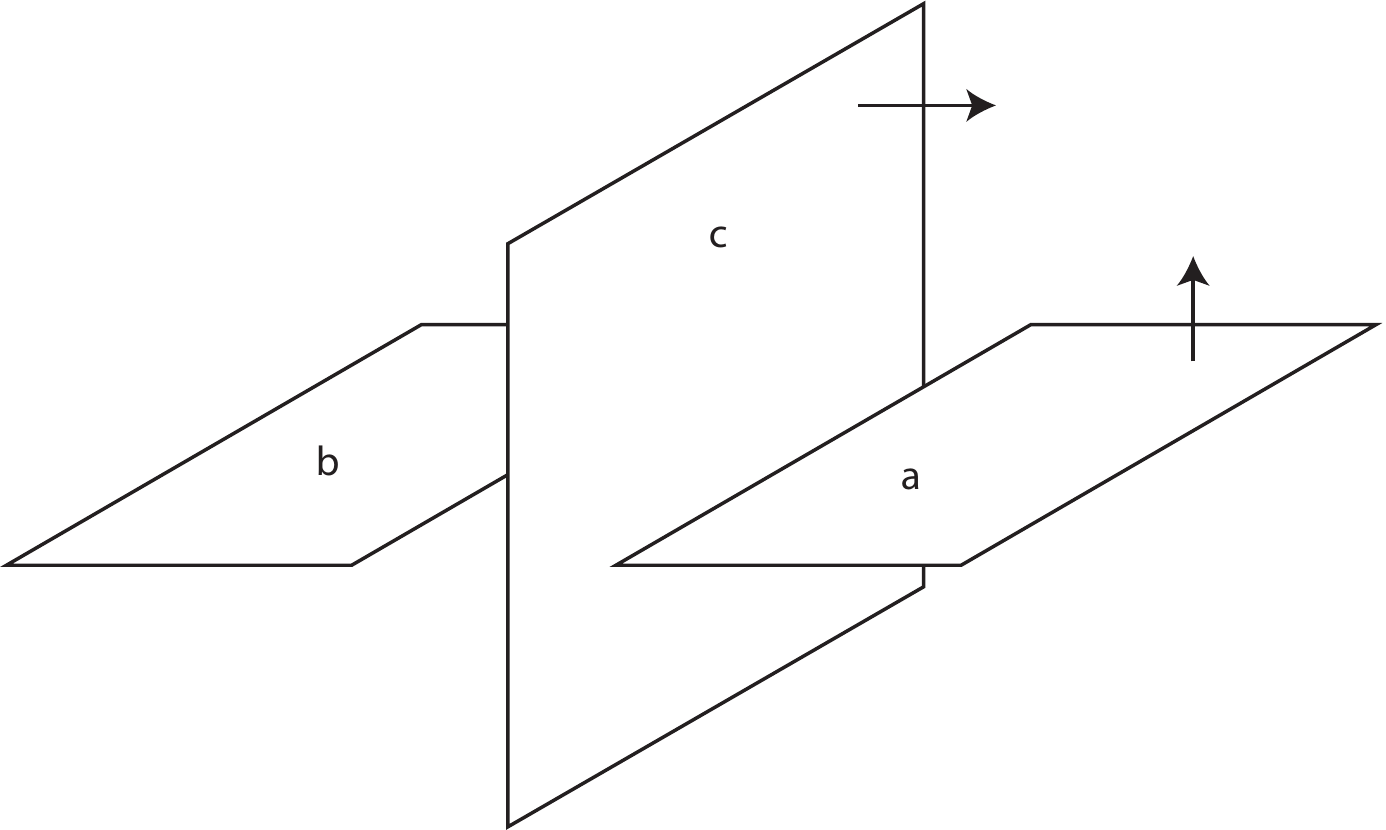}$ 
\caption{2-knot diagram leaves yielding a relation $a = b\circ c$.} \label{ris:quandle}
\end{center}
\end{figure}

Now for every double line we create a relation of the form  $a=b\circ c$ with $a,b,c \in Q$ corresponding to the connected components of the diagram as shown in Fig. \ref{ris:quandle}. Factorizing the free quandle modulo those relations one gets a quandle $Q_D$ for the diagram $D$. The following holds:

\begin{teo}
Consider a 2-knot $K$ and its diagrams $D_1, D_2$. The quandles $Q_{D_1}$ and $Q_{D_2}$ are isomorphic.
\end{teo}

In other words, the quandle is a knot invariant (for both 1- and 2-knots). This fact is verified by direct check of the moves.

The notion of parity allows one to refine the idea of quandle. To do that, one considers two operations --- $\circ$ and $\star$ --- and imposes certain conditions on them. After that all even singularities of diagram yield a relation of the form $a = b\circ c$ and the odd ones --- of the form $a=b\star c$.

We describe this strcture in more detail following the work of D.M. Afanasiev in 1-dimensional case (see \cite{afan}). First of all, let us define the notion of {\it parity quandle}.

\begin{opr}
A set $Q$ with two binary operations $\circ, \star: Q\times Q\to Q$ is called a {\rm parity quandle} if
\begin{enumerate}
    \item the operation $\circ$ is idempotent, i.e. $a\circ a = a$ for all $a\in Q$;
    \item the operations $\circ$ and $\star$ are right invertible, i.e. for all $b,c \in Q$ there exists a unique $a\in Q$ such that $$a\circ b = c$$ and there exists a unique $x\in Q$ such that $$x\star b = c;$$
    \item for all $a,b,c \in Q$ the following right distributivities hold:
    \begin{align} (a\circ b)\circ c = (a\circ c)\circ (b\circ c); \\ 
                  (a\circ b)\star c = (a\star c)\circ (b\star c); \\ 
                  (a\star b)\circ c = (a\circ c)\star (b\star c); \\
                  (a\star b)\star c = (a\star c)\star (b\circ c). 
    \end{align}
\end{enumerate}
\end{opr}

Note that parity quandle is not a quandle with respect to the operation $\star$. But if one takes $\star = \circ$ one gets a usual quandle.

Now for a given diagram $D$ of a 2-knot with parity we can construct a {\it parity quandle $Q_D$}. Again it is done via generators and relations.

As a set of generators we take the set of the diagram's leaves. Then every double line gives rise to a relation of the form $a\cdot b = c$ with $\cdot = \circ$ for even double lines and $\cdot = \star$ otherwise.

\begin{teo}
Parity quandle is an invariant of abstract 2-links; in other words equivalent diagrams have isomorphic parity quandles.
\end{teo}

It is important to note that quandle axioms (both usual and the parity one) are the same for all dimensions. \\

{\bf Example.} An example of the use of this invariant can also be obtained fro the variety of spun knots.

First, consider {\it decorated spun knots}: decorate every parallel with the sign ``+'' or ``--'', describing the neighbourhood of which parallel is ``upper'' and which is ``lower''. This decoration can be obtained from the decoration of the endpoints of the chords of the chord diagram yielding the spun knot.

Given a decorated spun knot one can calculate the corresponding parity quandle. Naturally this quandle is the same as the quandle for the underlying 1-knot. Thus two 1-knots with different parity quandles yield non-equivalent 2-knots.

\section{A list of problems}

One of the most important properties of parity for virtual 1-knots is the following.
There exists an invariant parity bracket $[K]$ for free links (hence for virtual links) valued in linear combinations of free knot diagrams. Diagrams in the bracket values are obtained from the diagram $K$ via {\it smoothings} and if the diagram $K$ is {\it odd and irreducible} (i.e. all its crossings are odd and one cannot apply a reducing second Reidemeister move to it) the following equality holds:

$$[K]=K.$$

Here on the left side $K$ is just a single diagram of a knot (i.e. equivalence class of diagrams) and on the right side --- as the only diagram with the coefficient 1.

That immediately leads to a number of corollaries. In particular if $K'$ is a diagram equivalent to the diagram $K$, then $[K']=K$.

That means that $K$ is obtained from $K'$ via smoothings. Thus not only we can state the nontriviality and minimality of the diagram $K$ but much more: in every equivalent diagram one can find the same picture.

This encourages one to seek similar structures in 2-dimensional case. Thus the following problems arise:

\begin{task}
Is an odd irreducible diagram of a 2-knot minimal?
\end{task}

\begin{task}
How one can construct picture-valued (diagram-valued) 2-knot invariants?
\end{task}

\begin{task}
To study group-valued parities.
\end{task}

\begin{task}
(L. Kauffman) Is the (sliced) genus of a virtual knot less or equal then its classical genus?
\end{task}

\begin{task}
Is every 2-knot with for which any parity is trivial classical?
\end{task}

Virtual knot theory can be applied to the study of classical knots. Historically the initial approach was the following. Let $K\sqcup L$ be a link consisting of a split link $K$ and a link $L$. Then $S^{3}\backslash K$ has a natural structure of Seifert fibration. Therefore the link $L$ can be projected to a Seifert surface. Knot diagrams on a 2-surface can have a non-trivial parity which can be used to create invariants of $K\sqcup L$.
 
 Thus virtual knot techniques, in particular the aprity theory, absent in classical case, may be applied to other (non-Reidemeister) theories of classical knots.
 
\begin{task}
It is well known that slice knots are very common in 2-knot theory. Thus a natural pronblem is to use similar structures to obtain non-trivial parities for 2-knot an link theories.	
\end{task}

\end{document}